\documentclass[a4paper,11pt]{article}
\usepackage{amsfonts,amssymb,latexsym,amsmath}
\usepackage{multicol}
\newcounter{num}[section]
\newcommand{\Num}{\refstepcounter{num}%
\textbf{\arabic{section}.\arabic{num}}}

\newcommand{\Theorem}{\textbf{Theorem~}}
\newcommand{\Proof}{\textbf{Proof.~}}
\newcommand{\Def}{\textbf{Definition}}

\newcommand{\Lemma}{\textbf{Lemma~}}
\newcommand{\Ex}{\textbf{Example~}}
\newcommand{\Rem}{\textbf{Remark~}}
\newcommand{\Prop}{\textbf{Proposition~}}
\newcommand{\Cor}{\textbf{Corollary~}}

\begin{document}

\Large

\newcommand{\gx}{{\mathfrak g}}
\newcommand{\nx}{{\mathfrak n}}
\newcommand{\hx}{{\mathfrak h}}

\newcommand{\ut}{{{\mathfrak u}{\mathfrak t}}}
\newcommand{\gl}{{{\mathfrak g}{\mathfrak l}}}
\newcommand{\UT}{{{\mathrm U}{\mathrm T}}}
\newcommand{\GL}{{{\mathrm G}{\mathrm L}}}

\newcommand{\Dp}{{\Delta^+}}
\newcommand{\De}{{\Delta}}
\newcommand{\Dt}{{\De^{[t]}}}
\newcommand{\tphi}{\tilde{\varphi}}
\newcommand{\vphi}{{\varphi}}

\newcommand{\al}{{\alpha}}
\newcommand{\La}{{\Lambda}}
\newcommand{\ga}{{\gamma}}
\newcommand{\Ga}{{\Gamma}}
\newcommand{\Xc}{{\cal X}}
\newcommand{\Oc}{{\cal O}}
\newcommand{\Ac}{{\cal A}}
\newcommand{\Bc}{{\cal B}}
\newcommand{\Ic}{{\cal I}}
\newcommand{\Jc}{{\cal J}}
\newcommand{\Yc}{{\cal Y}}
\newcommand{\Zc}{{\cal Z}}

\newcommand{\Mat}{{\mathrm{Mat}}}
\newcommand{\TTheta}{{\tilde{\Theta}}}
\newcommand{\TLa}{{\tilde{\La}}}
\newcommand{\TBc}{{\tilde{\Bc}}}

\newcommand{\Ab}{{\Bbb A}}
\newcommand{\Fb}{{\Bbb F}}
\newcommand{\Dc}{{\cal D}}
\newcommand{\wt}{w^{[t]}}
\newcommand{\wto}{w^{[t-1]}}

\newcommand{\row}{{\mathrm{row}}}
\newcommand{\col}{{\mathrm{col}}}
\newcommand{\Ad}{{\mathrm{Ad}}}
\newcommand{\VDp}{V_{D,\vphi}}

\renewcommand{\leq}{\leqslant}
\renewcommand{\geq}{\geqslant}

\title{Invariants of the coadjoint action on the basic varieties of the unitriangular group}
\author{A.N.Panov}

\date{}

\maketitle

\begin{abstract}
 We  find the generators of the fields of invariants of the coadjoint action of the unitriangular group on the basic varieties and basic cells. It is proved that the transcendental degree of the field of invariants on a basic cell  coincides with the number of factors  in the special factorization of the associated element of the Weyl group as a  product of reflections.
\end{abstract}

\section{Introduction}

For some finite   groups, the problem of classification of all irreducible representations is considered to be a "wild" problem \cite{Dr}. The unitriangular group over a finite field is one of these groups.
To construct the representation theory for such groups, the new program  arose in the paper \cite{DI}.
 In this program,   supercharacters and  superclasses play the some role as  irreducible characters and  conjugacy  classes in the usual representation theory. Emergence of this program was initiated by the theory of basic characters of the unitriangular group created by C.Andr\'e in a series of papers (see e.g.  \cite{A1,A2,A3}).

For constructing his theory, C.Andr\'e started from  A.A.Kirillov's orbit method for unipotent groups over a finite field $\Fb_q$. This method   establishes a  one to one correspondence
between the irreducible characters  and the coadjoint orbits (see \cite{K1}).
 For  a given basic subset $D\subset\Dp$ and  a map  $\vphi: D\to \Fb_q^*$, C.Andr\'e defined the basic variety $\VDp$  that is a subset  of the dual space  $\nx^*$ of the Lie algebra  $\nx =\ut(n,\Fb_q)$. The basic variety $\VDp$ is invariant with respect to the coadjoint representation of the group  $\UT(n,\Fb_q)$, and the dual space  $\nx^*$  splits into the basic varieties.  To each  basic variety $\VDp$, C.Andr\'e attached a basic character.
 This basic character
  decomposes into the irreducible characters  that correspond to the coadjoint orbits lying in the basic variety $\VDp$.

 In the paper \cite{A2}, C.Andr\'e introduced the notion of derived basic subset
       $D'$. To some roots from  $D'$, he  attaches  the invariant regular functions on the basic variety  $\VDp$.
       But these system of invariants is not complete in the sense
that it do not generate   the field of invariants.  Roughly speaking, it happens because using the  first derived subset  $D'$ one can define the second derived subset  $D''$ and the corresponding invariants, after, the third and so on  (see example \ref{second}).

In this paper, we consider basic varieties defined over a field of zero characteristic.  For a given basic subset  $D$, we construct its extension  $C(D)$.
For each  $\xi\in C(D)$, we determine the invariant rational function   $F_\xi$ on the basic cell $V^\circ_D$.
We prove that this  system of invariants is algebraically independent and generate the field of invariants on the basic cell (see Theorem  \ref{thm1}). Deleting the functions $F_\xi$, ~$\xi\in D$, in this system, we obtain  the system of generators on the basic variety  $\VDp $ (see Theorem  \ref{thm2}).

There are two approaches in construction of the subset $C(D)$.
 The  first approach uses  construction of diagram $\Dc$; with each root of $C(D)$, we associate the symbol "$\otimes$"\, on the diagram.  In the framework of the second approach, by subset $D$, we construct the  element (i.e permutation)  $w_D$ of the Weyl group (i.e. $S_n$). We characterize the elements of the form  $\{w_D\}$; the element $w\in W$ is an element of the form  $w_D$ if and only if the main affine neighborhood $\Yc_w$ of its Schubert variety   $\Xc_w$ is a cone with the center at origin
(see Theorem  \ref{hom} and Corollary  \ref{cor}).  We will refer to these elements $w\in W$ as homogeneous  elements. Using the permutation  $w_D$, we construct the subset $C(D)$ as follows;
 a root $\xi$ belongs to $C(D)$ if and only if the reflection $r_\xi$  appears in a certain factorization of  $w_D$  (see Theorem  \ref{dec}).

Notice that the main result of the paper \cite{pan}  is  a special case  of the one of this paper.

\section*{Homogeneous elements of the Weyl group}

Let  $K$ be a field of zero characteristic and
и $G=\GL(n,K)$. The unitriangular group  $N=\UT(n,K)$ is a subgroup of
$G$ that consists of all upper triangular matrices with ones on the diagonal.
Its Lie algebra  $\nx=\ut(n,K)$ consists of all upper triangular matrices with zeros on the diagonal.
Transposing  $\nx$, we obtain the subspace  $\nx_-$. Applying the invariant nondegenerate  form
$(A,B)=\mathrm{Tr}AB$ we identify  $\nx_-$ with the dual space  $\nx^*$ of $\nx$.

 Simplifying notations, we refer to any pair  $(i,j)$, where $1\leq
i,j\leq n$, ~ $i\ne j$, as a {\it root}. Here $i$ is the row number  of the root  $\al =
(i,j)$; respectively, $j$ is the column number. Denote
$i=\row(\al)$ and $j=\col(\al)$.

The Weyl group  $W=S_n$ acts on the set of all roots
корней $\De$ via the formula  $w(i,j) = (w(i),w(j))$.
We  refer to a pair  $(i,j)$, where
$n\geq i > j \geq 1$, as a {\it positive root}. We denote by $\Dp$ the subset of all positive roots.

We say that a root $\gamma=(i,k)$ is a sum of two roots $\al$ and $\beta$ if there exists $1\leq j\leq n-1$, ~ $j\ne i$,~ $j\ne k$ such that either $\al=(i,j)$, ~ $\beta=(j,k)$, or  $\beta=(i,j)$, ~ $\al=(j,k)$.

Following C.Andr\'e, we refer to a subset  $D\subset\Dp$  as a {\it basic} subset if there exists no more than one root of $D$ in any row and any column.

Let  $\gamma\in\Dp$. A root  $\al\in\Dp$ is
$\gamma$-{\it singular} if there exists a root  $\beta\in\Dp$ such that $\gamma=\al+\beta$.
 A positive root  $\al$ is $D$-{\it singular} if
there exists a positive root  $\beta$ such that  $\al+\beta\in
D$. We  denote by $S(D)$ the subset  of all  $D$-singular roots.

A root is $D$-{\it
regular} if it is not a  $D$-singular root. We denote the set of all $D$-regular roots by
$R(D)$. Then  $\Dp = S(D) \sqcup R(D)$.

By definition  $D\subset R(D)$. Denote  $M(D)=R(D)\setminus D$.
Notice that the set  $M(D)$ is closed with respect to addition  (i.e. if $\al,
\beta \in M(D)$ and  $\al+\beta $ is a root, then  $\al+\beta\in M(D)$).

With any basic subset   $D$ we associate  the element of the Weyl group  $w_D$. For each $1\leq j\leq n$,
the value  $w_D(j)$ equals to the greatest  $1\leq i \leq n$ that obeys the conditions\\
1)~ $(i,j)\notin M(D)$;\\
2)~ $i\notin\{w_D(1),\ldots, w_D(j-1)\}$.\\
{\bf Example}.\\
1) If  $D=\varnothing$, then  $S(D)=\varnothing$, ~$R(D)=M(D)=\Dp$, ~ $w_D=\mathrm{id}$.\\
2) Let $n=4$, ~$D=\{(3,1),~ (4,2)\}$.
Then \\$S(D)=\{(4,3),~ (3,2),~ (2,1)\}$, ~ $M(D)=\{(4,1)\}$, ~$w_D=\left(\begin{array}{cccc} 1&2&3&4\\
3&4&2&1\end{array}\right)$.\\
\Rem\Num\label{rem0}. Notice that, if  $(i,j)\in D$, then  $w_D(j)=i$. The converse statement is also true: if
$w_D(j)>j$, then  $(w_D(j),j)$ belongs to $D$.
 
 Let $\{E_{\ga}:~ \ga\in\Dp \}$ be the standard basis of matrix units in  $\nx_-$. For a given basic subset  $D$ and a map  $\vphi: D\to K^*$, we consider the element
\begin{equation}\label{xd} X_{D,\vphi} = \sum_{\xi\in D} \vphi(\xi) E_\xi\end{equation}
in $\nx_-=\nx^*$.

{\it A basic variety} $\VDp$ is the orbit of  $ X_{D,\vphi}$ with respect to the left-right action of the group  $N\times N$ on $\nx^*$ (see \cite{DI, Yan}).
Obviously, the basic varieties are invariant with respect the coadjoint representation  of the group   $N$.  The dual space  $\nx^*$ splits into basic varieties
   $$\nx^* = \bigsqcup \VDp.$$
Two basic varieties corresponding to the same  $D$ are  called {\it isotypical}.  We refer to the union of all isotypical basic varieties
$$V^\circ_D = \bigsqcup _{\vphi}\VDp $$
as a  {\it basic cell}. The basic cells form the partition of  $\nx^*$. {\it A basic cone} $V_D$ is the closure of a basic cell $V_D^\circ$ with respect to the Zariski topology.

 Let $B$ be the subgroup of upper triangular matrices in $\GL(n,K)$. Consider the flag variety  $\Xc= G/B$. The subset  $\Oc = N_-B\bmod B$  is called the {\it main affine neighborhood} of the point  $p=B\bmod B$. One can define the isomorphism  of affine varieties
$\Psi:\Oc\to\nx_-$ such that  $g = \left(1+\Psi(g)\right)\bmod B$
for any  $g\in \Oc$. The subspace  $\nx_-$ is naturally identified with  $K^{l}$, ~
where $l=\frac{n(n-1)}{2}$,  the point  $p$ is identified with the origin of  $K^{l}$.

The flag variety  $\Xc$ splits into the Schubert cells
$$\Xc_w^\circ =  B\dot{w}B\bmod B,~\mbox{~~ where ~~}w\in W.$$
{\it The Schubert variety} $\Xc_w$ is the closure of the Schubert cell
$\Xc_w^\circ$ with respect to the Zariski topology. It is well known that every Schubert variety contains the origin  $p$. For any   $w\in W$, we introduce the notations:
$$ \Yc_w = \Psi(\Xc_w\cap \Oc),\quad \quad\Yc^\circ_w = \Psi(\Xc_w^\circ\cap \Oc).$$
It is easy to see that  $\Yc_w$ is the closure of  $\Yc_w^\circ$.\\
\Def. An element  $w\in W$ is  {\it homogeneous}, if
$\Yc_w$ is a cone  in  $\Oc$ with the center at the point $p$.\\
\Theorem\Num\label{hom} The element  $w\in W$ is homogeneous if and only if $w=w_D$ for some basic  subset
$D$.\\
\Proof
For every   $n\times n $-matrix $A$,    $w\in W$ and
$1\leq i,j\leq n$, we will use notations:
\\
i)~ $J'=\{1\leq b\leq j-1 :~ w(b)>i\}$,~~ $J = J'\cup\{j\}$;\\
ii)~ $I'=w(J')$,~~ $I = I'\cup\{i\}$; \\
iii)~ $P_{ij}(A)$ is a minor of the matrix $A$  with the system of rows $I$ and columns $J$; if $\eta=(i,j)$ is a root (i.e.
 $i\ne j$) we will denote  this minor by  $P_\eta$;\\
iv)~ $I_{w,j}$ is the subset of all  $w(j)\leq k \leq n$ that $k\ne w(b)$ for each  $1\leq b\leq j-1$.

It is not difficult to prove that a matrix $A$ belongs to the Schubert cell  $\Xc^\circ_w$ whenever  $P_{w(j),j}(A)\ne 0$ for  each
$1\leq j\leq n$, and
$$ P_{ij}(A) = 0~~\mbox{for~ each}~~ i\in I_{w,j}.$$
Therefore, a matrix  $X\in \nx_-$ belongs in  $\Yc_w^\circ$ whenever $P_{w(j),j}(1+X)\ne 0 $ for each $1\leq j\leq n$,
 and
\begin{equation}\label{ex} P_{ij}(1+X) = 0~~\mbox{for~ each}~~ i\in
I_{w,j}.\end{equation} 1)~ Let  $ w=w_D$ for some  basic subset  $D$. Then  $i\in I_{w,j}$ implies  $(i,j)\in M(D)$.
Hence, $i>j$ and the equality
(\ref{ex}) transforms to  \begin{equation} P_{ij}(X) = 0~~\mbox{for ~
each}~~ i\in I_{w,j}.\end{equation}  Therefore,  $\Yc_w$ is a cone.\\
2) Let  $w\ne w_D$ for any basic subset  $D$. Then there exists a pair  $(i,j)$ such that  $i\leq j$ and $i\in I_{w,j}$.
Assume that  here  $j$ is the least number and  $i$ is the greatest number among all pairs with this property.

As we said above, the polynomial  $P=P_{ij}$ is zero on  $1+\Yc_w$.
Suppose that  $|I'|=m$. Then  $|J'|=m$ and $|I|=|J|=m+1$. Hence, for each  $X\in \nx_-$, we obtain
$$P(1+X) = \sum_{s=0}^{m+1} P^{(s)}(X),$$
where $ P^{(s)}(X)$ is a homogeneous polynomial of degree  $s$.

Assume that $\Yc_w$ is  a cone. Then, all its homogeneous components   $ P^{(s)}(X)$ are zero on  $\Yc_w$.

There is natural projection  $\pi: \Mat(n,K)\to \Mat(m+1,K)$ (respectively,   $\pi': \Mat(n,K)\to \Mat(m,K)$) that sends every  $n\times n$-matrix  to its submatrix with the system of rows  $I$ and the system of columns  $J$ (respectively, $I'$ и $J'$). Then,  for every matrix  $A\in \Mat(n,K)$, one can obtain  its image  $\pi'(A)$ by removing the first row and the last column of the matrix  $\pi(A)$. Notice also that  $P(A)=\det \pi(A)$.

Denote $k=j-i+1$. The image  $\pi(N_-)$ consists of all  $(m+1)\times(m+1)$-matrices of the form
\begin{equation}\label{im}
Y =\left(\begin{array}{cc} Y_{11}&Y_{12}\\
Y_{21}&Y_{22}
\end{array}\right),
\end{equation}
where $Y_{12}$ is a lower triangular matrix  of size  $k\times k$. Similarly,  the image
$\pi'(N_-)$ consists of all $m\times m$-matrices of the form   (\ref{im}) with  $k'\times k'$-block  $Y_{12}$, where $k'=k-2$ for $k
\geq 2$, and  $k'=0$ for $k=1$.

The image $\pi'(\dot{w})$ is a monomial matrix in  $\GL(m,K)$. Then,  $\pi'(B\dot{w}B)$
is a Bruhat class of the element  $\pi'(\dot{w})$ in $\GL(m,K)$. It is well known that the closure of each  Bruhat class contains  the Borel subgroup. Therefore,  the closure of the Bruhat class  $\pi'(B\dot{w}B)$ contains the Borel subgroup of  $\GL(m,K)$. Hence, the closure of  $\pi'(1+\Yc_w)$ contains the subset  of all  $m\times m$ upper triangular matrices  $Y'$ of the form (\ref{im}) with the corresponding  $k'\times k'$-block $Y_{12}$. The closure of the image  $\pi(1+\Yc_w)$ contains the subset  $\Yc$
that consists of all  $(m+1)\times(m+1)$-matrices  $Y$ obeying the conditions:\\
1) $Y$ has the form  (\ref{im}) with corresponding $k\times k$-block  $Y_{12}$,\\
2) removing the first row and the last column of the matrix  $Y$, we obtain an  upper triangular matrix,\\
3) $P(Y)=0.$

According to assumption,  $ P^{(s)}(\Yc_w) = 0$ for each  $0\leq s\leq m+1$.
Then $ P^{(s)}(\Yc) = 0$. Recall that $P(Y)=\det(Y)$. It is easy to prove that the homogeneous components of $\det(Y)$ are not zero on $\Yc$.  Therefore,  $\Yc_w$ is not a cone.
$\Box$\\
\Cor\Num \label{cor}  The basic cone  $V_D$ coincides with the cone  $\Yc_{w_D}$. Basic cones are the tangent cones at origin of the Schubert varieties of the homogeneous elements of the Weyl group.
\\
 \Proof The generators  of defining ideal  of the basic cell  $V_D^\circ$ (see \cite{A1}) coincide with the ones of  the Schubert cell  $\Yc_{w_D}^\circ$ (see the proof above). This proves the statement. $\Box$

\section*{Extension of the basic subsets and associated diagrams}

Fix the linear order $\succ$ on the set of all positive roots as follows:~
~$(k,m)\succ(i,j)$, if $m<j$ or $m=j$,~ $k>i$. According to this order
\begin{equation}\label{order}(n,1)\succ(n-1,1)\succ\ldots\succ(2,1)\succ(n,2)\succ\ldots \succ(n,n-1).
\end{equation}
For a given  basic subset $D$, we construct the diagram  $\Dc$ that is a $n\times n$-matrix;  the places over and on the diagonal  are not filled,
and the places below the diagonal  are filled with the symbols  "$\otimes$"{}, "$\bullet$"{}, "$+$"\, и "$-$"\ according to some rules listed below. We associate with each  positive root  $(i,j)$, where
$i>j$,  the place  $(i,j)$ on the diagram  $\Dc$.
The places  $(i,j)\in M(D)$ are filled with the symbol "$\bullet$"{}. The places filled with the symbol "$\otimes$"\, form the subset
$$C(D)=\{\xi_1\succ \xi_2\succ\ldots\succ\xi_c\}\subset \Dp.$$
We start  to construct  the diagram by filling all places  of
$M(D)$ with the symbol "$\bullet$". We refer to this procedure as a zero step of  diagram construction.

 We put the symbol "$\otimes$"\, on the greatest (in the sense of  $\succ$) place  $\xi_1$ from $\Dp\setminus M(D)$.
   If $\xi_1=\al+\beta$, where  $\col(\al) = \col(\xi_1)$ and $\row(\beta) = \row(\xi_1)$,
then we fill the place  $\al$ with the symbol  "$+$"\,, and the place  $\beta$  with  "$-$".
 This step  finishes the first in diagram construction.

Further, we put the symbol  "$\otimes$"{} on the greatest  (in the sense of order  $\succ$) empty place from  $\Dp$.
Denote the corresponding positive root by $\xi_2$. Let  $\xi_2 =
\al + \beta$, where $\col(\al) = \col(\xi_2)$ and $\row(\beta) = \row(\xi_2)$.

 If the places  $\al$ and $\beta$ are empty after the previous steps, then we fill the place  $\al$ with the symbol  "$+$"\, and the place $\beta$ with "$-$"\,. If one of these two steps is already filled, then we do not fill the other place.

 After filling with  "$+$"\ and  "$-$"\,, we finish the next step завершается calling it as the second step. Continuing the procedure,
 we obtain the diagram $\Dc$.  The number of the last step coincides with the number of symbols "$\otimes$"{} on the diagram. The roots corresponding the symbols "$\otimes$"\, form the subset  $C(D)$.

For each  $\gamma\in \Dp$, we denote by  $A_\gamma$ the set of all
$\gamma'\in C(D)$ such that  $\row(\ga')=\row(\ga)$ and $1\leq \col(\ga')
<\col(\ga)$. In the following Lemma all roots are positive.\\
\Lemma\Num\label{lem1}
a) Let  $\xi=\xi_i\in C(D)$, ~$\xi=\al+\beta$, where
$\col(\al)=\col(\xi)$ and $\row(\beta)=\row(\xi)$. Suppose that  the places  $\al$ and $\beta$
were empty on the diagram before the   $i$th step. Then $A_\al\ne\varnothing$, if  $A_\xi\ne\varnothing$.\\
b) The root  $\xi$ belongs to  $D$ if and only if  $\xi\in C(D)$ and
$A_\xi = \varnothing$. \\
c) Let  $\xi=\xi_i\in C(D)$, ~$\xi=\al+\beta$, where
$\col(\al)=\col(\xi)$ and $\row(\beta)=\row(\xi)$. Then, either one of the places  $\al$ and $\beta$ is filled with the symbol  "$+$"\, or "$-$"\, before the  $i$th step and the other place is empty, or both places are empty before the  $i$th step.
\\
\Proof
We will prove this three statements simultaneously using the induction method with respect to the number of  $\xi$
in the list of positive roots (\ref{order}). For the greatest root, the statements are obvious. Assume that the statements are true  for all roots  $\succ \xi$; we will prove it for  $\xi$.\\
\textbf{a)} Suppose that  $\xi$, $\al$, $\beta$ satisfy conditions of the statement a).
 Then, at the  $i$th step,  $\xi$ is filled with "$\otimes$"\,,    $\al$ (respectively, $\beta$) is filled with "$+$"\, (respectively,  "$-$"\,).

Assume that   $A_\xi\ne\varnothing$ and $A_\al=\varnothing$.  Since  $A_\xi\ne\varnothing$, there exists
$\xi_{i_0}\in C(D)$,~ $i_0<i$, lying on the left side and in the same row as  $\xi$.  Decompose
$\xi_{i_0}=\al'+\beta$, where $\row(\al')=\row(\al)$ (see the table below). Since  the place
$\beta$ is filled with "$-$"\, at the  $i$th step  and  $i_0<i$, the place  $\beta$ is  empty after  the  $i_0$th step.
At the same $i_0$th step, the place  $\al'$ can't be empty (in this case,  $\al'$ is filled with "$+$"\, and
$\beta$ with "$-$"\, at the   $i_0$th step). The place  $\al'$ can't be filled with "$\otimes$"\, (obvious),
can't be filled with "$-$"\, (as $A_\al\ne\varnothing$),
and can't be filled with "$\bullet$"\, (this contradicts to the induction assumption of the statement c) applied for $\xi_{i_0}$).

 Therefore,  there exists  $\xi_{i_1}=\al'+\beta'\in C(D)$,
where $i_1<i_0<i$,  lying lower and in the same column as  $\xi_{i_0}$. At the $i_1$th step, the place  $\al'$ is filled with
the symbol  "$+$"\, and the place  $\beta'$ with "$-$"\,. Then,  $\xi_{i_1}=\xi_0 + \beta''$, where
$\row(\beta'') = \row(\xi_{i_1})$. Since the place $\xi_{i_0}$
is filled with "$\otimes$"\, at the  $i_0$th step, it is empty at the $i_1$th step.
Therefore,  the place  $\beta''$ is filled  (with  "$-$") at some step  $i_2<i_1$. There exists
$\xi_{i_2}\in C(D)$ lying on the left side and in the same row as $\xi_{i_1}$.  For reader's convenience, we show all  these roots on a table
 linking two roots by the horizontal segment if two roots lie in the same row, and by
 vertical segment if two roots lie in the same column.

$$\begin{array}{ccccccccc}&&\al'&\rule{0.5cm}{0.1mm}&\al&&&&\\
&&\rule{0.1mm}{0.5cm}&&\rule{0.1mm}{0.5cm}&&&&\\
&&\xi_{i_0}&\rule{0.5cm}{0.1mm}&\xi&\rule{0.5cm}{0.1mm}&\beta&&\\
&&\rule{0.1mm}{0.5cm}&&&&\rule{0.1mm}{0.5cm}&&\\
\xi_{i_2}&\rule{0.5cm}{0.1mm}&\xi_{i_1}&&&&\beta'&\rule{0.5cm}{0.1mm}&\beta''
\end{array}$$

Notice that decomposition  $\xi_{i_1}=\al'+\beta'$ obeys the condition a)
of Lemma. For the root  $\xi_{i_1}$, there exists  an element of
$C(D)$ lying  on the left side and in the same row (this is  $\xi_{i_2}$), but  for  $\al'$ there are no such elements. This contradicts the induction assumption for  the statement  a).\\
\textbf{b)} Let us prove the statement  b).\\
{\bf b1}) Let  $\xi\in D$.  Suppose that  $A_\xi\ne \varnothing$.
According to the induction assumption,  the leftmost root of  $A_\xi$ belongs to  $D$.
But there can't be two roots of $D$ lying in the same row. Therefore $A_\xi=\varnothing$.

Let us show that  $\xi\in C(D)$. Assume that  $\xi\notin C(D)$.
Hence, the place  $\xi$ can't be filled on the diagram  $\Dc$ with   "$\otimes$"\,.
Also it can't be filled with "$\bullet$"\, (as $\xi\in D$), can't be filled with "$-$"\, (as
 $A_\xi=\varnothing$).
Suppose that it is filled with "$+$"\,. There exists a root   $\xi_{j}\in C(D)$ such that   $\xi_{j} = \xi
+\beta$;  at the $j$th step, $\xi$ is filled with "$+$"\, and $\beta$
  with  "$-$"\,. The root  $\xi_{j}$ don't belong to  $D$
(otherwise two roots of  $D$ will lie in the same column).  The induction assumption for statement b) implies
 $A_{\xi_{j}}\ne \varnothing$.  From the induction assumption for statement a), we obtain
 $A_\xi\ne \varnothing$. This contradicts what was proved above.\\
{\bf b2}) Let   $\xi = \xi_i\in C(D)$ and $A_\xi = \varnothing$. Let us prove that  $\xi\in D$. Assume that  $\xi\notin D$. Since
$A_\xi = \varnothing$,  according to the induction assumption, there are no  roots of $D$ on the left side  and in the same row as  $\xi$. If there is no roots of  $D$  lower and in the same column as $\xi$, then  $\xi\in M(D)$. This contradicts  the condition   $\xi\in C(D)$.
 Therefore, there is  root  $\ga\in D$ lying lower and in the same column as  $\xi$.
Then, according to induction assumption, $\gamma=\xi_j\in C(D)$,~ $j<i$, and $A_\ga=\varnothing$. Then, at the $j$th step,   the place  $\xi$ is filled with "$+$"\,; this contradicts
$\xi\in C(D)$. Finally, $\xi\in D$.
 \\
\textbf{c)} Let  $\xi=\xi_i$, ~$\al$,~ $\beta$ satisfy  the condition of the statement c).
Suppose that the symbols    "$+$"\, and  "$-$"\,  don't appear on the places $\al$ and $\beta$ before the  $i$th step. Let us show that, in this case, these places are empty  before the
  $i$th step  (then at the $i$th
 $\al$ is filled with "$+$"\, and  $\beta$ is filled with "$-$"\,).

 Indeed,  the symbol "$\otimes$"\,  can't appear on the places  $\al$ and  $\beta$ before the  $i$th step. Since  $\xi\in C(D)$, either  $\xi\in D$, or there exists a root of $D$ on the left side and in the same row  (according induction  assumption, it will the leftmost root from
  $A_\xi$).
 All these implies that there is no symbol "$\bullet$"\, on the place  $\beta$.

  Suppose that  the place  $\al$ is filled with "$\bullet$"\,  and the place $\beta$ is empty.
  Then  $\xi\notin D$ ($D$-singular roots can't lie in  $M(D)$).
  The item  b) implies that  $A_\xi\ne\varnothing$. Further, arguing similarly as in item  a), one can show that  $A_\al\ne\varnothing$. The leftmost root of  $A_\al$ belongs to  $D$; this leads to a contradiction
  (the symbols  "$\bullet$"\, and "$\otimes$"\, can't lie in the same row). $\Box$\\
\Prop\Num\label{prop1}
1) Each root of  $D$ belongs to  $C(D)$. \\
2) The basic subset  $D$ consists of those roots of  $\xi\in C(D)$ that  $A_\xi\ne\varnothing$.\\
3) Each root of $D$ is sited  in his column over all roots from   $C(D)$.
\\
\Proof
The items  1) and 2) are corollaries of the statement  b) of the previous Lemma.
Suppose that the root   $\xi\in C(D)$ lies over some root   $\eta\in D$. Then   $\eta=\xi+\beta$, ~ $\col(\eta)=\col(\xi)$,~ $\row(\eta)=\row(\beta)$. Since
$A_\eta =\varnothing$, the place  $\beta$ is empty at the step corresponding to $\eta$. Then, at this step  $\xi$ is filled with the symbol  "$+$"\, and  $\beta$ with "$-$". This contradicts to  $\xi\in C(D)$. $\Box$

\section*{Decomposition of homogeneous elements of the Weyl group}

As in the previous section, we extend  a basic subset  $D$ to the subset $C(D)=\{\xi_1\succ\ldots\succ\xi_c\}$.
The subset $C(D)$ consists of all roots that are filled on the diagram  $\Dc$ with the symbol
"$\otimes$"\,. For each  $1\leq i\leq c$, denote
$w_i = r_1\cdots r_i$ where  $r_j$ is the reflection with respect to the root $\xi_j$. Set
$w_0=e$.
\Prop\Num\label{prop2} Let $\xi_i\in C(D)$ and  $\col(\xi_i)=t$.
Let  $\ga\in\Dp$ and $\col(\ga)\geq t$. After the $i$th step, the place  $\ga$ is empty or  is filled with "$\bullet$"\,  if and only if
$w_i(\ga)>0$. Respectively, after the $i$th step, the place  $\ga$  is filled with one of the symbols  "$\otimes$"\,, "$+$"\,,"$-$"\,  if and only if  $w_i(\ga)<0$.
\\
\Proof We apply the induction method with respect to  $i$. For  $i=0$ the statement is obvious.
 Assume that the statement is true for numbers $< n$. Let us prove for  $n$.

Simplifying notations, we set $\xi=\xi_i$. Notice that, by induction assumption,  $w_{i-1}(\xi)>0$.

We prove the statement in each of the following five cases separately:\\ 1) $\xi\succ\gamma$ and $(\xi,\ga)>0$,~ 2) $\ga\succ \xi$ and
$(\xi,\ga)>0$,~ 3)~ $(\xi,\ga)<0$,~ 4) ~$\ga=\xi$, ~
5)~$(\ga,\xi)=0$.\\
{\bf 1)}~ $\xi\succ\gamma$ and $(\xi,\ga)>0$. In this case,  $\xi$ is a sum of two positive roots $\xi=\ga+(\xi-\ga)$. Then
\begin{equation}\label{w1}
w_i(\ga)=w_{i-1}r_i(\ga) = -w_{i-1}(\xi-\ga).
\end{equation}

Suppose that the place  $\ga$  is empty after the  $i$th step or  is filled with "$\bullet$"\,.
According to the statement  c)  of Lemma  \ref{lem1}, the place $ \xi-\ga$ is filled before the $i$th step with  "$+$"\, or "$-$"\,. According to the induction assumption, we have  $w_{i-1}(\xi-\ga)<0$. The formula (\ref{w1}) implies $w_i(\ga)>0$.

Suppose that $\ga$ is filled with one of the symbols  "$+$"\, or "$-$"\, after the  $i$th step.
We have two cases. The first case is when this symbol appears at the $i$th step, the second case is when it appears before the $i$th step.

In the first case, the both places  $\ga$ and $\xi-\ga$ are empty before the  $i$th step. By induction assumption, $w_{i-1}(\ga)>0$
and  $w_{i-1}(\xi-\ga)>0$. The formula  (\ref{w1})
implies  $w_i(\ga)<0$.

In the second case, according to the induction assumption,  $w_{i-1}(\ga)<0$. Recall that    $w_{i-1}(\xi)>0$.
 The equality   $\xi = \ga+(\xi-\ga)$ implies  $w_{i-1}(\xi-\ga)>0$. By (\ref{w1}), we obtain   $w_i(\ga)<0$.\\
{\bf 2)}~ $\ga\succ \xi$ and $(\xi,\ga)>0$. In this case, $\ga$ is a sum of two positive roots $\ga=\xi+(\ga-\xi)$, where
$\col(\ga)=\col(\xi)$ and $\row(\ga) = \row(\ga-\xi)$.  Then
\begin{equation}\label{w2}
w_i(\ga)=w_{i-1}r_i(\ga) = w_{i-1}(\ga-\xi).
\end{equation}
Since $\ga\succ \xi$, the place  $\ga$ can't be empty before the $i$th step. If the place  $\ga$ is filled with "$\bullet$"\,, then there is no "$\otimes$"\, on the left side and in the same row.   Then, before the $i$th step,  the place  $\ga-\xi$
is empty or is filled with "$\bullet$"\,. According to the induction assumption,  $w_{i-1}(\ga-\xi)>0$.
The formula  (\ref{w2}) implies $w_i(\ga)>0$.

Suppose the  $\ga$ is filled with one of the following symbols
"$\otimes$"\,, "$+$"\, or "$-$"\, after the $i$th step. Since  $\ga\succ \xi$,  we filled $\ga$ before the $i$th step.
By induction assumption, $w_{i-1}(\ga)<0$. The equalities  $\ga=\xi+(\ga-\xi)$ and $w_{i-1}(\xi)>0$ imply
$w_{i-1}(\ga-\xi)<0$. By the formula (\ref{w2}), we conclude $w_i(\ga)<0$.\\
 {\bf 3)}~ $(\xi,\ga)<0$. In this case,  $\xi+\ga$ is a root and  $\col(\xi+\ga)=\col(\xi)$, ~  $\row(\xi+\ga) = \row(\ga)$. Then
\begin{equation}\label{w3}
w_i(\ga)=w_{i-1}r_i(\ga)=w_{i-1}(\xi+\ga).
\end{equation}

Suppose that $\ga$ is empty or is filled with "$\bullet$"\, after the  $i$th step.
In this case, $\ga$ is empty or is filled with "$\bullet$"\, after the  $(i-1)$th step.
According the induction assumption,  $w_{i-1}(\ga)>0$.
Since $w_{i-1}(\xi)>0$, we obtain   $w_i(\ga)=w_{i-1}(\xi+\ga) = w_{i-1}(\xi)+
w_{i-1}(\ga)>0$.

As $\col(\ga)>\col(\xi)$, the place $\ga$ can't be filled with "$\otimes$"\, or "$+$"\,. Suppose that
$\ga$ is filled with the symbol  "$-$"\, after the  $i$th step. Since
$\row(\ga)>\row(\xi)$, the place $\ga$ is filled before the  $i$th step. The place lies
$\xi+\ga$ lower and in the same column as  $\xi$. Therefore, the place  $\xi+\ga$ is also filled  before the  $i$th step.
It can't be filled with "$\bullet$"\, (in this case,  $\ga$ can't be filled with "$-$"\,). Hence $w_{i-1}(\xi+\ga)<0$.
Therefore $w_i(\ga)<0$.\\
{\bf 4)} ~$\ga=\xi$. In this case, $\ga$ is filled with "$\otimes$"\, and
$$w_i(\ga)=w_{i-1}r_i(\xi)=-w_{i-1}(\xi)<0.$$
{\bf 5)}~ $(\ga,\xi)=0$. We have $r_i(\ga)=\ga$, the statement follows from the induction assumption. $\Box$

Recall that $ w_i = r_1 \ldots r_i$ and $C(D)=\{\xi_1\succ \xi_2\succ\ldots\succ\xi_c\}$.\\
\Theorem\Num\label{dec} The permutation $w_D$ decomposes into a product of reflections
$ w_D = r_1 \ldots r_c$ where $r_i=r_{\xi_i}$. In this decomposition, \\
i) ~ $w_i(\xi_j)>0$ for all $1\leq i<j$,\\
ii) ~ $ \xi_j\notin M(D)$,\\
iii) each $\xi_j$ is the greatest root
(in the sense of  $\succ$) that obeys i), ii) and is less than $\xi_{j-1}$.\\
\Proof  We have to prove that  $w_c=w_D$, i.e.   $w_c(t)=w_D(t)$ for each $1\leq
t\leq n$. For $t=1$, the statement is obvious. Assume that
$w_c(j)=w_D(j)$ for each  $1\leq j < t$. Let us prove the statement for   $j=t$.

Denote by  $\wt$ the product of all  $r_\xi$, where $\xi\in
C(D)$ and $\col(\xi)\leq t$; the factors  are ordered as in  $C(D)$. The permutation  $\wt$ coincides with $w_a$,
where $\xi_a$ is the least root in  $\{\xi\in C(D):~ \col(\xi)\leq t\}$.

Notice that  $\wt(j)=w_c(j)$ for each  $1\leq j\leq t$  (as  $r_\eta(j)=j$ for each  $\eta\in\Dp$, ~ $\col(\eta)>t$).
According to the induction assumption,  $\wt(j)=w_D(j)$ for each  $1\leq j
< t$. Let us prove that  $\wt(t)=w_D(t)$. Denote
$$ A = \{\wt(m):~ t\leq m\leq n\} = \{w_D(m):~ t\leq m\leq n\},$$
$$ M = \{m: ~t< m \leq n, ~ (m,t)\in M(D)\}.$$

The  both numbers  $\wt(t)$ and $w_D(t)$ belong to $A$, and  $w_D(t)$,
by definition, is the greatest number in  $A\setminus M$.

We have to show that   $\wt(t)$ is the greatest number in  $A\setminus M$.
We present  $A$ as a union of three disjoint subsets
\begin{equation} \label{adec} A = \{\wt(t)\}\bigsqcup \{\wt(m): m\in M\}\bigsqcup \{\wt(m): m\notin M\}.
\end{equation}
If  $m\in M$, then the root  $\ga=(m,t)$, lies in  $M(D)$.
Then there is no elements of $C(D)$  in the  $m$th row  and on the left side with respect to $\ga$. Hence $\wt(m)=m$.
We see that the second subset in decomposition (\ref{adec}) coincides with  $M$.

Suppose that  $t< m \leq n$ and $m\notin M$. Then  $\ga=(m,t)\notin
M(D)$. After filling the $t$th column, the place  $\ga$ is filled with one of the symbols "$\otimes$"\,,
"$+$"\, or "$-$"\,. The Proposition \ref{prop2}  implies
$\wt(\ga)<0$. Therefore,  $ \wt(m) < \wt(t)$. Then,  $\wt(t)$ is the greatest number of  $A\setminus M$. That is $\wt(t)=w_D(t)$. $\Box$

\section*{Invariants of the basic varieties}

Recall some  definitions of the theory of Poisson algebras.
\emph{A Poisson algebra} is a commutative associative algebra with unit equipped  with a Poisson bracket. We refer to an ideal  $\Ic$ of a Poisson algebra   $\Ac$ as   \emph{a Poisson ideal} if $\{\Ic,\Ac\}\subset \Ic$. An element  $a\in \Ac$ is called  {\it a Casimir element } if $\{a,\Ac\}=0$.
{\it A tensor product} of two Poisson algebras  $\Bc_1$ и $\Bc_2$ is the algebra  $\Bc_1\otimes \Bc_2$ equipped with a Poisson bracket such that  $\{\Bc_1,\Bc_2\} = 0$, and such that $\Bc_1$, $\Bc_2$ are  Poisson subalgebras.

{\it The standard Poisson algebra } is a Poisson algebra   $\Ab_n$ generated by  $$p_1, \ldots, p_n, q_1,\ldots, q_n$$ obeying the relations  $\{p_i,q_j\}=\delta_{ij}$. Properties of Poisson algebra  $\Ab_n$ are similar to properties of the Weyl algebra
    (see subsection  4.6 of the book  \cite{Dix}). In particular, it is known that  every derivation of the Poisson algebra   $\Ab_n$ is inner.  Every Poisson ideal of the Poisson algebra  $\Ab_n\otimes \Bc$ has the form  $\Ab_n\otimes \Jc$, where $\Jc$
     is a Poisson ideal of  $\Bc$. Every Casimir element of  $\Ab_n\otimes \Bc$ belongs to  $\Bc$.\\
     \Lemma\Num \label{lem12} Let  $\Ac$ be an arbitrary Poisson algebra, $p$,~ $q$ are two element of  $\Ac$ such that  $\{p,q\}=1$.
Let  $\Ac^\diamond  = \{ a\in \Ac: ~ \{a,q\}=0\}$. Suppose that \\
 1) the derivation  $D_p(a) = \{p,a\}$  is locally nilpotent;\\
 2) $p$ and $\Ac^\diamond$ generate the   $\Ac$ as a commutative associative algebra.\\
 Then, the Poisson algebra  $\Ac$ is isomorphic to  $\Ab_1\otimes \Bc$, where $\Bc=\Ac^\diamond/\Ac^\diamond q$.
 \\
 \Proof[Proof scheme] The proof is similar to Lemma  4.7.5 of \cite{Dix}.
 Consider $D_p:\Ac^\diamond \to\Ac^\diamond $. Define a map $\Theta_p:\Ac^\diamond \to\Ac^\diamond$ by the formula
\begin{equation}\label{di}
\Theta_p(a) = \sum_{s=0}^\infty (-1)^sD^s_p(a)\frac{q^s}{s!}.
\end{equation}
The map  $\Theta_p$ is a homomorphism of the Poisson algebra $\Ac^\diamond$.
In addition, $\{p,\Theta_p(a)\} = \{q,\Theta_p(a)\} = 0$ for every  $a\in \Ac^\diamond $.
The kernel of  $\Theta_p$ coincides with
$\Ac^\diamond q$.  It is easy to prove that  $\Ac= \Ab_1\otimes \mathrm{Im} \Theta_p$. $\Box$

  Let $\gx$ be an arbitrary Lie algebra, $G$ be its adjoint group. The algebra   $K[\gx^*] $  is a Poisson algebra with respect to the Poisson bracket induced by the commutator of  $\gx$.  In the case of  $\gx=\gl(n,K)$, the Poisson bracket has the form
            $\{x_{ij}, x_{km}\} = \delta_{jk}x_{im}$. Obviously, an ideal  $\Ic$  in $K[\gx^*] $ is a Poisson ideal if and only if
             the ideal $\Ic$ is invariant with respect to the coadjoint representation of the group $G$. Respectively, a
             Casimir element in  $K[\gx^*] $  is an invariant of the  coadjoint representation.

Let  $\xi=\xi_i=(s,t)$.  Consider two subsets of the system of positive roots:
$$\La_i =\{\eta\in \Dp:~ \col(\eta)\geq t, ~  w_{i-1} (\eta) > 0\},$$
$$\La_{>i} =\{\eta\in \Dp:~ \col(\eta)\geq t, ~  w_{i} (\eta) > 0\}.$$

That is the  root  $\eta$  belongs to  $\La_i$ (respectively,  $\La_{>i}$), if $\col(\eta)\geq t $ and the place  $\eta$ is empty on the diagram  or is filled with "$\bullet$"\, before (respectively, after) the $i$th step (see Proposition  \ref{prop2}). Notice that $\La_{>i}\subset \La_{i}$, and the both subsets are subalgebras in $\Dp$.

The linear span  $\nx_i$ of the system  $\{E_{-\eta},~~\eta\in \La_i\}$ is a Lie subalgebra of $\nx$. One can identify the dual space  $\nx_i^*$ with  the set  $\{\sum x_\eta E_{\eta}\}$, where $\eta$ runs through  $\La_i$. The algebra of regular functions  $\Ac_i=K[\nx^*_i]$ on  $\nx^*_i$ is a Poisson algebra; it can be identified with a polynomial algebra  $K[x_\eta: \eta\in \La_i]$. Similarly, we define the Lie algebra  $\nx_{>i}$ and the Poisson algebra  $\Ac_{>i}=K[\nx^*_{>i}]$.
The Lie algebra   $\nx_{>i}$ is a subalgebra of  $\nx_i$; using projection  $\nx_i^*\to\nx_{>i}^*$, one can identify  $\Ac_{>i}$
with the subalgebra of  $\Ac_i$.

Consider the ideal  $\Jc_i$ in $\Ac_i$ generated by the elements $x_\eta$, where  $\eta\in\La_i$ and $\row(\eta)> s$.
The construction method of $C(D)$ implies that if  $\eta\in\La_i$ and $\row(\eta)> s$, then  $\eta\in M(D)$.
It is easy to see that  $\Jc_i$ is an Poisson ideal in  $\Ac_i$.

Similarly, we define  the Poisson ideal  $\Jc_{>i}$ of $\Ac_{>i}$ as an ideal generated by all  $x_\eta$, where  $\eta\in\La_{>i}$ and $\row(\eta)> s$.
 Denote by  $\Bc_i$ a factor algebra of $\Ac_i$ with respect to the ideal $\Jc_i$; and  by $\Bc_{>i}$  a factor algebra of $A_{>i}$ with respect to $\Jc_{>i}$. Notice that $x_\xi$ is a Casimir element of  $\Bc_i$.

 Let  $\Bc_i'$
be a localization of  $\Bc_i$ with respect to the denominator subset generated by  $x_\xi$.
Decompose the set  $$T=\{\beta\in \La_t: ~ \col(\beta)>\col(\xi),~ \row(\beta)=\row(\xi)\}$$ into subsets $T=T_0 \cup T_+ \cup T_-$ where
$$ T_0=\{\beta\in T:~ \beta\in \La_t,~ \xi-\beta\notin \La_t\},$$
$$ T_+=\{\beta\in T:~ \beta\in \La_t,~ \xi-\beta\in \La_t\},$$
$$ T_-=\{\beta\in T:~ \beta\notin \La_t\}.$$
By the way, the places of  $T_-$ are filled with "$-$"\, after all steps before filling   $\xi$.
Denote  $ S_\pm = \{\al\in \Dp:~~ \xi-\al\in T_\pm\}$. If $$T_+\cup T_- = \{\beta_k\succ\ldots\succ\beta_1\}~ ~\mbox{and}
~~ \al_i = \xi-\beta_i,$$ then $S_+\cup S_- = \{\al_1\succ\ldots\succ\al_k\}$.

For the elements   $p_i=x_{\beta_i}$ and $q_j=x_{\al_j}x_\xi^{-1}$, we have $\{p_i,q_j\}=\delta_{ij}$.
Denote by  $\Ab_+$ (respectively,  $\Ab_-$) the standard Poisson algebra  generated by  $p_i,~q_i$ where $\al_i\in S_+$ and $\beta_i\in T_+$ (respectively, $\al_i\in S_-$ and $\beta_i\in T_-$).\\
\Lemma\Num\label{lem2}
1) The Poisson algebra  $\Bc_i'$ is isomorphic to the Poisson algebra   $$\Ab_+\otimes K[x^{\pm 1}_{\xi}]\otimes \Bc_{>i}.$$
2) The isomorphism  between these algebras defines a Poisson embedding  \begin{equation}\label{tts}\Theta_{i}: \Bc_{>i}\hookrightarrow  \Bc_i'~~ \mbox{obeying}~~ \{\Theta(x),\Ab_+\}=0.\end{equation}
 In addition, for every  $\eta \in\La_{>i}$,  the image of  $x_\eta$ has  the form
\begin{equation}\label{ttt}
\Theta_{i}(x_\eta) = x_\eta + Q_{\succ \eta}\bmod \Jc_i,\end{equation}
where  $Q_{\succ \eta}$ is a polynomial in  $x_\ga$,~ $\ga\in \La_i$, ~ $\ga\succ \eta$, and $x_{\xi}^{-1}$.

Notice that  $\{\Theta(x),x_\xi\}=0$, since  $x_\xi$ is a Casimir element in $\Bc_i'$.
\\
\Proof
It is sufficient to prove that there exists a Poisson embedding
\begin{equation}\label{Theta}
   \Theta: \Bc_{>i}\hookrightarrow \Bc_i'
\end{equation}
obeying  (\ref{tts}) and (\ref{ttt}).

 Let us prove that  $\La_i\cup T_-$ is a subalgebra in $\Dp$.
Really, let  $\eta_1\in\La_i$,~ $\eta_2\in T_-$ and $\eta_1 + \eta_2\in\Dp$. There are two cases.

The first case: $\eta_2 = (s,a)$, ~ $ \eta_1 = (a,b)$, where $s>a>b$. Then  $\eta_1 + \eta_2 =(t,b)\in T \subset
\La_i\cup T_-$.

 The second case: $\eta_1 = (b,s)$, ~ $ \eta_2 = (s,a)$, where $b>s>a$. Then   $\eta_1 + \eta_2 =(b,a)$ and  $\xi + \eta_1=(b,t)\in\La_i$.
Suppose that  $\eta_1+\eta_2\notin \La_i$. Then,  the place  $(b,a)$ is filled with "$-$"\, at the step before filling
     $\xi$. The place $(b,t)$ is filled at that step with the symbol  "$\bullet$"\, (otherwise,
      $(b,t)$ is filled with "$\otimes$"\,,  $\xi$ with "$+$"\,, and $\eta_1$ with "$-$"\,). Since  $(b,a)$ is filled with "$-$"\, and  $t<a$, there is  "$\otimes$"\, on the left side  and in the same row as  $(b,t)$. Then, more left there is some "$\otimes$"\, on the place belonging to $D$. This contradicts to definition of the subset $M(D)$. So  $\eta_1 + \eta_2\in\La_i$.

The linear span  $\tilde{\nx}_i$ of the system  $\{E_\eta:~\eta\in \La_i\cup T_-\}$ is a Lie subalgebra of  $\nx$.
The algebra of regular functions  $\tilde{\Ac}_i$ on $\tilde{\nx}_i^*$   is a Poisson algebra.
The Lie algebra  $\nx_i$ is a subalgebra of   $\tilde{\nx}_i$; using the natural projection   $\tilde{\nx}_i^*\to\nx_i^*$, one can identify $\Ac_i$ with the subalgebra of $\tilde{\Ac}_i$.
Algebra $\tilde{\Bc}_i$ and  its localization  $\tilde{\Bc}_i'$ are defined similarly as  $\Bc_i$ и $\Bc_i'$.

Denote by  $\TLa_{>i}$ the subset  $\La_{>i}\setminus S_-$. Let
$\TBc_{>t}$ be a subalgebra of $\Bc_{>i}$ generated by  $x_\eta$, where $\eta\in \TLa_{>i}$.

It is not difficult to prove  that the pair $p_1,q_1$ of  $\tilde{\Bc}_{i}'$
obeys conditions of Lemma  \ref{lem12}. Then, the pair  $p_2,q_2$ also obeys the same conditions in the factor algebra of
  $(\tilde{\Bc}_{i}')^\diamond$ modulo the ideal generated by  $q_1$ and so on. Applying  Lemma \ref{lem12} for pairs $(p_1,q_1), \ldots,(p_k,q_k)$ , we obtain
  \begin{equation}\label{tild}
\tilde{\Bc}_{i}' \cong \Ab_+\otimes \Ab_-\otimes K[x^{\pm 1}_\xi]\otimes \tilde{\Bc}_{>i}.
    \end{equation}
 The isomorphism (\ref{tild}) defines a Poisson embedding
  $\tilde{\Theta}:\tilde{\Bc}_{>t}\hookrightarrow \tilde{\Bc}_t'$ obeying  (\ref{tts}) and (\ref{ttt}). Therefore,
 \begin{equation}\label{brat}
\TTheta(x_{\eta_1+\eta_2}) =
\{\TTheta(x_{\eta_1}),\TTheta(x_{\eta_2})\},
\end{equation}
for all $\eta_1$ and $\eta_2$ from  $\TLa_{>t}$.

By direct calculations, we obtain formulas for  $\TTheta(x_\eta)$.
If  $\eta = (a,b)$, where $a<s$ and $b>t$, then $\xi=\beta+\eta+\al$, where $\beta=(s,a)$ and $\al = (b,t)$.
Then
\begin{equation}\label{theta1}
  \TTheta(x_\eta) = - x_\xi^{-1}\cdot\left|\begin{array}{cc} x_{\al + \eta}&x_\eta\\
  x_{\xi}&x_{\eta+\beta}\end{array}\right| = x_\eta - pq,
        \end{equation}
where $p=x_{\eta+\beta}$ and  $q=x_{\al+\beta}x_\xi^{-1}$.
If $\eta = (a,s)$, where $a>s$, then
\begin{equation}\label{theta2}
    \TTheta(x_\eta) =  x_\eta + x_\xi^{-1}\cdot \sum_{\beta\in T_+\cup T_-} x_{\beta+\eta}x_{\xi-\beta}
    \end{equation}
In all other cases  $\TTheta(x_\eta) =x_\eta$.

Let $\eta\in\La_{>i}$.
We say that  $\eta$ is a root of the first type if  $\eta\in \TLa_{>t}$  and the image  $\TTheta(x_{\eta})$  belongs to  ${\Bc}_i'$. Respectively,  $\eta$ is a root of the second type if $\eta\in S_-$,  or $\eta\in \TLa_{>t}$ and the image $\TTheta(x_{\eta})$
 don't belong to  ${\Bc}_i'$.
 Easy to see that  the  root $\eta $ is a root of the second type if and only if $\eta$ belongs to  $S_-$ or   $\eta=(a,b)$ as in formula  (\ref{theta1}) with  $\eta+\beta\in T_-$.
Define the map  $\Theta$ as follows:
$$\Theta(x_\eta) = \left\{\begin{array}{cccccc} \TTheta(x_\eta), & \mbox{if}&~ \eta~ &~ \mbox{is~ a~ root~of~the~} &~ \mbox{first}~ &~ \mbox{type},\\
x_\eta, & \mbox{if}&~ \eta~ & \mbox{is~ a~ root~ of~ the~} &~ \mbox{second}~ &~ \mbox{type}.
\end{array}\right. $$
Let us show that    \begin{equation}\label{bra}\Theta(x_{\eta_1+\eta_2}) =
\{\Theta(x_{\eta_1}), \Theta(x_{\eta_2})\}.\end{equation}

Really, if $\eta_1$,~ $\eta_2$ are the roots of the  first type, then  $\eta_1+\eta_2$ is a root of the first type, and   (\ref{bra}) follows from  (\ref{brat}). If $\eta_1$,~ $\eta_2$ are the roots of the second type, then   $\eta_1+\eta_2$ is also a root of the second type; formula
 (\ref{bra}) is obvious. If one of roots, say  $\eta_1$, has the first type  and  the other, i.e.
     $\eta_2$, has the second type, then it is easy to prove that  $\eta_1+\eta_2$ has the first type.
Let   $\eta_2=(a,b)$ be as in the formula  (\ref{theta1}).
The formulas    $\{\TTheta(x_{\eta_1}), p\}=\{\TTheta(x_{\eta_1}), q\}=0$  imply
  $$\Theta(x_{\eta_1+\eta_2}) = \TTheta(x_{\eta_1+\eta_2}) =
\{\TTheta(x_{\eta_1}),\TTheta(x_{\eta_2})\} =$$ $$ \{\Theta(x_{\eta_1}), x_{\eta_2}-pq\}=
\{\Theta(x_{\eta_1}), x_{\eta_2}\} =\{\Theta(x_{\eta_1}), \Theta(x_{\eta_2})\}.$$
One can check that the conditions 2) holds for  $\Theta$.
$\Box$

Let  $\Ac^\circ$ be a localization of the algebra of regular functions  $\Ac=K[\nx^*]$ with respect to the denominator subset
 generated by the minors  $P_\xi$, где $\xi\in D$ (see section 2).
Consider the ideal  $\Ic_D^\circ$ of the algebra $\Ac^\circ$ generated by $P_\ga$, where $\ga\in M(D)$. Then \begin{equation}\label{vd}K[V^\circ_D] = \Ac^\circ/\Ic_D^\circ.\end{equation}
The ideal  $\Ic=\Ic^\circ\cap\Ac$ is the   defining ideal of the basic cone   $V_D$. Hence  $K[V_D] = \Ac/\Ic_D$.\\
\Theorem\Num\label{thm1}. The field of invariants of the coadjoint action of the group $N$ on the basic  cell  $V^\circ_D$ (respectively, basic cone $V_D$)  is a pure transcendental extension of the main field  $K$ of degree  $|C(D)|$.
One can choose a basis  $\{F_\xi:~ \xi\in C(D)\}$ of the field of invariants   such that, if  $\xi\in D$, then  $F_\xi$ coincides with the minor $P_\xi$.\\     
\Proof
We will construct the  elements  $\{F_\xi\}$ by the induction process with respect to the descending order $\succ$.
 Suppose that we have   defined  the elements  $F_{\xi_1},\ldots, F_{\xi_{i-1}}$.
We  choose    $F_{\xi_1}=x_{\xi_1}\in \Ac$;  the  element  $F_{\xi_2}$ belongs to the algebra  $\Ac^{(1)}$ that is the localization of $\Ac$ with respect to $ F_{\xi_1}$; the element $F_{\xi_3}$ belongs to the algebra $\Ac^{(2)}$ that is the localization of  $\Ac^{(1)}$ with respect to  $ F_{\xi_2}$ and so on.
Denote by  $\Ac^{(i-1)}$ the algebra that appears as a result of consistent localization of  $\Ac$ with respect to $F_{\xi_1},\ldots, F_{\xi_{i-1}}$.

For a given  $\eta\in\Dp$,  let  $\xi_{i-1}$ be the greatest (in the sense of  $\succ$) root such that  $\xi_{i-1}\succ \eta$.
Denote by  $F_\eta $ an arbitrary element of the subset
$$ \Theta_1\circ\cdots\circ\Theta_{i-1}(x_\eta).$$
The element   $F_{\eta}$ is defined  uniquely modulo the ideal  $\Jc^{(i-1)}$ generated by all $F_\eta$, where $\eta\in M(D)$, $\eta\succ\xi_{i-1}$. One can present  $F_{\eta}$ in the form
 \begin{equation}\label{zzz} F_{\eta} = x_{\eta}+ Q_{\succ\eta},\end{equation}
 where  $Q_{\succ\eta}$  is a polynomial in  $x_\ga$, ~$\ga\succ \eta$, and $F_{\xi_1}^{-1},\ldots, F_{\xi_{i-1}}^{-1}$ (see Lemma \ref{lem2}).   Localizing the algebra  $\Ac^{(i-1)}$ with respect to $F_{\xi_i}$, we obtain the algebra  $\Ac^{(i)}$
and its ideal  $\Jc^{(i)}$. Applying  \ref{lem2},  we have
\begin{equation}\label{at}
\Ac^{(i)}/\Jc^{(i)} \cong \Ab_{m_i}\otimes K[F^{\pm 1}_{\xi_1},\ldots,F^{\pm 1}_{\xi_i}]\otimes \Bc_{>i}.
\end{equation}

Continuing the process, we obtain the algebra  $\Ac'= \Ac^{(c)}$ and its ideal  $\Jc'$ (the other notation is  $\Jc^{(c+1)}$) generated by all  $F_\eta$, where  $\eta\in M(D)$. Then,
\begin{equation}\label{last}
\Ac'/\Jc' \cong \Ab_{m}\otimes K[F^{\pm 1}_\xi: ~ \xi\in C(D)].
\end{equation}

The formula  (\ref{zzz}) implies that\\
1) ~ $\Jc'$  is a prime ideal of  $\Ac'$; respectively,  $\Jc=\Jc'\cap\Ac$
is a prime ideal of  $\Ac$;\\
2) ~ the field of rational Casimir functions in  $\mathrm{Frac}(\Ac/\Jc)$ coincides with the field of rational functions in
$F_{\xi_i}$ where $1\leq i\leq c$.

 We have to prove that  $\Jc=\Ic$.  Set  $\Ic^{(0)}=\Jc^{(0)} =\{0\}$. Let $\Ic^{(i)}$ be the ideal of $\Ac^{(i)}$ generated by  $P_\eta$, where $\eta\succ \xi_i$. Let   $\Ic'$ (further also  $ \Ic^{(c+1)}$) be  the localization of the ideal  $\Ic$ with respect to  $F_{\xi_1}, \ldots, F_{\xi_c}$. Using the induction method on  $0\leq i\leq c$, we will prove that \\
i) $\Ic^{(i)} =  \Jc^{(i)}$;\\
ii) if $\xi=\xi_i\in D$, then  $P_{\xi} = F_{\xi} M_{>\xi}\bmod \Ic^{(i-1)}$, where $M_{>\xi}$ is a product of some
$P_{\xi'}^{\pm}$, ~$\xi'\in D$,~ $\xi'\succ \xi$.

 For $i=0$ the statements are trivial. Assume that the statement are true for $i$, let us prove for  $i+1$.

 Let $\xi_{i}\succ \eta_1\succ\ldots\succ\eta_m\succ \xi_{i+1}$ be the system of all elements of  $M(D)$ lying between  $\xi_i$
 and $\xi_{i+1}$.
  Let $\Ic^{(i,k)}$ (respectively,  $\Jc^{(i,k)}$) be the ideal of  $\Ac^{(i)}$ generated by  $\Ic^{(i)}$ (respectively, $\Jc^{(i)}$) and $P_{\eta_j}$ (respectively, $F_{\eta_j}$), where $1\leq j\leq k$.
Set  $\Ic^{(i,0)} = \Ic^{(i)}$ (respectively, $\Jc^{(i,0)}=\Jc^{(i}$). To prove i), we will show that  $\Ic^{(i,k)} = \Jc^{(i,k)}$ for each  $1\leq k\leq m$. Apply the induction method on $k$.  For $k=0$ this is obvious.
Assume that the statement is proved for  $k$; let us prove for  $k+1$.

 Let  $\eta=\eta_{k+1}$. Expanse the the minor  $P_\eta$ along the last column:
 \begin{equation}\label{ppp}
 P_\eta = x_\eta A_\eta+Q_{>\eta},
  \end{equation}
   where $A_\eta$ is the cofactor of  $x_\eta$, and $Q{>\eta}$ has the form as in  (\ref{zzz}). Easy to prove that   $A_\eta$ is invertible modulo $\Ic$ in the localization of the algebra  $\Ac$ with respect to  $P_{\xi'}$, ~ $\xi'\in D$ ,~ $\xi'\succ \eta$.
    Using the  induction assumption for statement ii), we conclude that  $A_\eta$ is invertible in  $\Ac^{(i-1)}$.
Denote  $P'_\eta=\Ac_\eta^{-1}P_\eta$; the formula  (\ref{ppp}) implies
$P'_\eta = x_\eta + Q'_{>\eta}$, where $Q'_{>\eta}$ as in (\ref{zzz}).
To prove  ii), it is sufficient  to show that  $P'_\eta = F_\eta$.
Really, let $R_{>\eta} = P'_\eta - F_\eta = Q'_{>\eta} - Q_{>\eta}$.
The elements  $P'_\eta$  and $F_\eta$  are  Casimir elements modulo the ideal  $\Ic^{(i,k)} =  \Jc^{(i,k)}$.
 Then, $R_{>\eta}$ is also a  Casimir element; it can be expressed by the elements  $x_\ga$, ~ $\ga\succ\eta$.
By the formula   (\ref{at}), we obtain
\begin{equation}\label{att}
\Ac^{(i)}/\Ic^{(i,k)} \cong \Ab_{m_i}\otimes K[F^{\pm 1}_{\xi_1},\ldots,F^{\pm 1}_{\xi_i}]\otimes \left(\Bc_{>i}/\sim\right),
\end{equation}
where  $\Bc_{>i}/\sim$ is a factor algebra of  $\Bc_{>i}$ with respect of ideal generated by  $x_{\eta_1}, \ldots,x_{\eta_k}$.

Applying the formula (\ref{att}) and properties of  $R_{>\eta}$, we have
\begin{equation}\label{rrr}R_{>\eta}\in K[F^{\pm 1}_{\xi_1},\ldots,F^{\pm 1}_{\xi_{i}}]\bmod \Ic^{(i)}.\end{equation}
The elements  $P'_\eta$ and $F_\eta$ have weight  $\eta$ with respect to the action of the Cartan subalgebra $\hx$.
Then,  $R_{>\eta}$ has the same weight  $\eta$. Let $\eta=(a,b)$. Consider the element  $h_s$ of the standard basis of  $\hx$
 that has one on the  $a$th place  and zeroes on the other places. Then,  $\mathrm{ad}_{h_a} R_{>\eta} = - R_{>\eta} $.
On the other hand,  since  $\eta\in M(D)$, there is no element  $\xi'\in C(D)$ lying in the $a$th row and on the left side from $\eta$. In addition, there is no  element  $\xi'\in C(D)$ that  $\xi'\succ \eta$ and $\col(\xi')=a$. Hence, $\mathrm{ad}_{h_a} x_{\xi'} = 0$ for every  $\xi'\succ\eta$, ~$\xi'\in C(D)$. By the formula (\ref{rrr}), we conclude that  $\mathrm{ad}_{h_a}R_{>\eta}=0$. This proves that  $R_{>\eta}=0$ and the equality $\Ic^{(i,k+1)} = \Jc^{(i,k+1)}$.
Taking  $k=m$, we prove i).

The statement  ii) can be proved similarly to  i). Then,  taking  $i=c+1$, we obtain  $\Ic'=\Jc'$. Finally,  $\Ic=\Jc$.

The statement  ii) enables to substitute  $F_\xi$, ~$ \xi\in D$, for $P_\xi$ in the list of generators. $\Box$\\
\Theorem\Num\label{thm2}.  The field of invariants of the coadjoint action of the group $N$ on the basic  variety
$V_{D,\phi}$   is a pure transcendental extension of the main filed  $K$ of degree   $|C(D)|-|D|$.\\
\Proof
 The statement is a corollary of previous theorem.

\section*{Examples}

\Ex\Num \label{first} Consider  $n=8$ and $D=\{(4,1), ~(7,2), ~ (8,3), ~
(5,4)\}$
 We construct the diagram  $\Dc$ in five steps (beginning from the zero step).

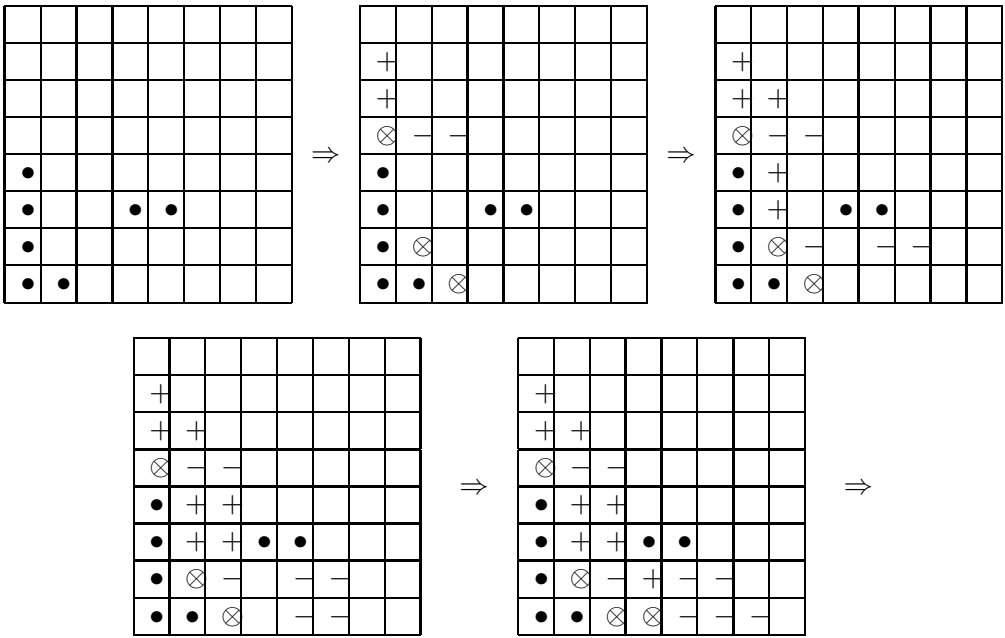
\begin{center}
\normalsize {
\begin{tabular}{|p{0.05cm}|p{0.05cm}|p{0.05cm}|p{0.05cm}|p{0.05cm}|p{0.05cm}|p{0.05cm}|p{0.05cm}|}
\hline & & &  & & &  & \\
\hline & & & &  & &   &\\
\hline & & & &  & &  &\\
\hline  & & & & & &  &\\
\hline $\bullet$ &  &   &   & &  && \\
\hline $\bullet$ &  & &$\bullet$ &$\bullet$ &  &&\\
\hline  $\bullet$ &   &  &  & & &&\\
\hline   $\bullet$&  $\bullet$ &   &  &  & && \\
\hline
\end{tabular}
 ~$\Rightarrow$~
\begin{tabular}{|p{0.05cm}|p{0.05cm}|p{0.05cm}|p{0.05cm}|p{0.05cm}|p{0.05cm}|p{0.05cm}|p{0.05cm}|}
\hline & & &  & & &  & \\
\hline $+$& & & &  & &   &\\
\hline $+$& & & &  & &  &\\
\hline $\otimes$ &$-$ &$-$ & & & &  &\\
\hline $\bullet$ &  &   &   & &  && \\
\hline $\bullet$ &  & &$\bullet$ &$\bullet$ &  &&\\
\hline  $\bullet$ & $\otimes$  &  &  & & &&\\
\hline   $\bullet$&  $\bullet$ &$\otimes$   &  &  & && \\
\hline
\end{tabular}
 ~$\Rightarrow$~
\begin{tabular}{|p{0.05cm}|p{0.05cm}|p{0.05cm}|p{0.05cm}|p{0.05cm}|p{0.05cm}|p{0.05cm}|p{0.05cm}|}
\hline & & &  & & &  & \\
\hline $+$& & & &  & &   &\\
\hline $+$&$+$ & & &  & &  &\\
\hline $\otimes$ &$-$ &$-$ & & & &  &\\
\hline $\bullet$ & $+$ &   &   & &  && \\
\hline $\bullet$ & $+$ & &$\bullet$ &$\bullet$ &  &&\\
\hline  $\bullet$ & $\otimes$  &$-$  &  &$-$ &$-$ &&\\
\hline   $\bullet$&  $\bullet$ &$\otimes$   &  &  & && \\
\hline
\end{tabular}}
\end{center}

\begin{center}
\normalsize {
\begin{tabular}{|p{0.05cm}|p{0.05cm}|p{0.05cm}|p{0.05cm}|p{0.05cm}|p{0.05cm}|p{0.05cm}|p{0.05cm}|}
\hline & & &  & & &  & \\
\hline $+$& & & &  & &   &\\
\hline $+$&$+$ & & &  & &  &\\
\hline $\otimes$ &$-$ &$-$ & & & &  &\\
\hline $\bullet$ & $+$ & $+$  &   & &  && \\
\hline $\bullet$ & $+$ &$+$ &$\bullet$ &$\bullet$ &  &&\\
\hline  $\bullet$ & $\otimes$  &$-$  &  &$-$ &$-$ &&\\
\hline   $\bullet$&  $\bullet$ &$\otimes$   &  &$-$  &$-$ && \\
\hline
\end{tabular}
\quad $\Rightarrow$\quad
\begin{tabular}{|p{0.05cm}|p{0.05cm}|p{0.05cm}|p{0.05cm}|p{0.05cm}|p{0.05cm}|p{0.05cm}|p{0.05cm}|}
\hline & & &  & & &  & \\
\hline $+$& & & &  & &   &\\
\hline $+$&$+$ & & &  & &  &\\
\hline $\otimes$ &$-$ &$-$ & & & &  &\\
\hline $\bullet$ & $+$ & $+$  &   & &  && \\
\hline $\bullet$ & $+$ &$+$ &$\bullet$ &$\bullet$ &  &&\\
\hline  $\bullet$ & $\otimes$  &$-$  &$+$  &$-$ &$-$ &&\\
\hline   $\bullet$&  $\bullet$ &$\otimes$   &$\otimes$  &$-$  &$-$ &$-$& \\
\hline
\end{tabular}}
 \quad $\Rightarrow$ \quad
 \end{center}
 \vspace{0.3cm}
\begin{center}
\normalsize {  $\Dc$ =
\begin{tabular}{|p{0.05cm}|p{0.05cm}|p{0.05cm}|p{0.05cm}|p{0.05cm}|p{0.05cm}|p{0.05cm}|p{0.05cm}|}
\hline & & &  & & &  & \\
\hline $+$& & & &  & &   &\\
\hline $+$&$+$ & & &  & &  &\\
\hline $\otimes$ &$-$ &$-$ & & & &  &\\
\hline $\bullet$ & $+$ & $+$  & $\otimes$  & &  && \\
\hline $\bullet$ & $+$ &$+$ &$\bullet$ &$\bullet$ &  &&\\
\hline  $\bullet$ & $\otimes$  &$-$  &$+$  &$-$ &$-$ &&\\
\hline   $\bullet$&  $\bullet$ &$\otimes$   &$\otimes$  &$-$  &$-$ &$-$& \\
\hline
\end{tabular}}
\end{center}
\vspace{0.2cm}
The subset  $C(D)$ consists of  $D$ and the root  $(8,4)$.
In this example,  $\{(8,4)\}$ coincides with the derived subset  $D'$ in the sense of Andr\'e (see \cite{A2})
and  $C(D) = D\cup D'$.
To each root of $C(D)$, we attach the polynomial
$F_{4,1} = x_{41}$,~  $F_{7,2} = x_{72}$,~  $F_{8,3} = x_{83}$,~
 $F_{5,4} = \left|\begin{array}{ccc} x_{52}&x_{53}&x_{54}\\
 x_{72}&x_{73}&x_{74}\\
 x_{82}&x_{83}&x_{84}\end{array}\right |$, ~ $F_{8,4} = x_{84}x_{41}+x_{83}x_{31}$.

The basic sell $V_D^\circ$ consists of all $X\in \nx_-$ obeying relations
$$ x_{51}=x_{61}=x_{71}=x_{81}=x_{82}=0,\quad \left|\begin{array}{ccc} x_{62}&x_{63}&x_{64}\\
 x_{72}&x_{73}&x_{74}\\
 x_{82}&x_{83}&x_{84}\end{array}\right | =\left|\begin{array}{ccc} x_{62}&x_{63}&x_{65}\\
 x_{72}&x_{73}&x_{75}\\
 x_{82}&x_{83}&x_{85}\end{array}\right | = 0,
$$
$$ F_{4,1} \ne 0 ,~  F_{7,2} \ne 0,~  F_{8,3} \ne 0,~  F_{5,4} \ne 0.$$
The field of invariants  of the basic cell  $V_D^\circ$
is the field of rational functions of $F_{4,1},~  F_{7,2},~ F_{8,3}, ~F_{5,4},~ F_{8,4}$.

Let $F_{i,j}^\circ$ be a value of $F_{i,j}$ at $X_{D,\varphi}$ (see (\ref{xd})).
The basic variety $V_{D,\varphi}$ consists of all $X\in V_D^\circ$ obeying relations
$F_\xi = F_\xi^\circ$, where $\xi\in D$.
The field of invariants  of the basic variety
 $V_{D,\vphi}$ is the field  of rational functions of  $F_{8,4}$.
 \\
 \Ex\Num \label{second} Let  $n=7$  and $D=\{(4,1),~ (5,2), ~(6,3),~ (7,5)\}$.
We construct the diagram  $\Dc$ in six steps (beginning from the zero step).
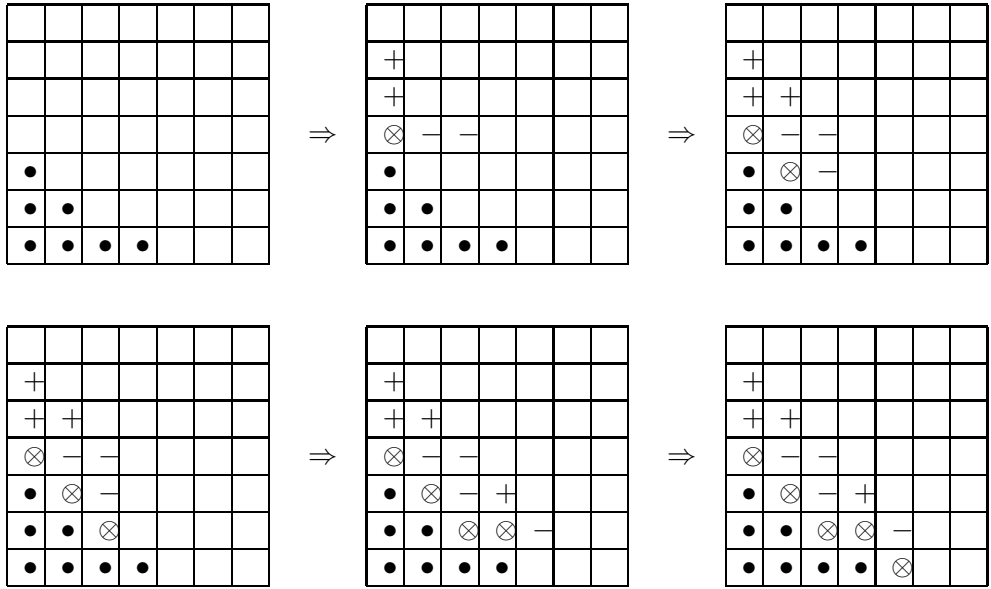
\begin{center}
\normalsize {
\begin{tabular}{|p{0.07cm}|p{0.07cm}|p{0.07cm}|p{0.07cm}|p{0.07cm}|p{0.07cm}|p{0.07cm}|}
\hline & & &  & & &   \\
\hline & & & &  & &  \\
\hline & & & &  & &  \\
\hline  & & & & & &  \\
\hline $\bullet$ &     &   & &  && \\
\hline $\bullet$ & $\bullet$  & & &  &&\\
\hline  $\bullet$ & $\bullet$  &$\bullet$  &$\bullet$  & & &\\
 \hline
\end{tabular}
\quad $\Rightarrow$\quad
\begin{tabular}{|p{0.07cm}|p{0.07cm}|p{0.07cm}|p{0.07cm}|p{0.07cm}|p{0.07cm}|p{0.07cm}|}
\hline & & &  & & &   \\
\hline $+$& & & &  & &  \\
\hline $+$& & & &  & &  \\
\hline $\otimes$ &$-$ &$-$ & & & &  \\
\hline $\bullet$ &     &   & &  && \\
\hline $\bullet$ & $\bullet$  & & &  &&\\
\hline  $\bullet$ & $\bullet$  &$\bullet$  &$\bullet$  & & &\\
 \hline
\end{tabular}
\quad $\Rightarrow$\quad
 \begin{tabular}{|p{0.07cm}|p{0.07cm}|p{0.07cm}|p{0.07cm}|p{0.07cm}|p{0.07cm}|p{0.07cm}|}
\hline & & &  & & &   \\
\hline $+$& & & &  & &  \\
\hline $+$&$+$ & & &  & &  \\
\hline $\otimes$ &$-$ &$-$ & & & &  \\
\hline $\bullet$ & $\otimes$    & $-$  & &  && \\
\hline $\bullet$ & $\bullet$  & & &  &&\\
\hline  $\bullet$ & $\bullet$  &$\bullet$  &$\bullet$  & & &\\
 \hline

\end{tabular}}
\end{center}
\vspace{0.3cm}
\begin{center}
\normalsize {
\begin{tabular}{|p{0.07cm}|p{0.07cm}|p{0.07cm}|p{0.07cm}|p{0.07cm}|p{0.07cm}|p{0.07cm}|}
\hline & & &  & & &   \\
\hline $+$& & & &  & &  \\
\hline $+$&$+$ & & &  & &  \\
\hline $\otimes$ &$-$ &$-$ & & & &  \\
\hline $\bullet$ & $\otimes$    & $-$  & &  && \\
\hline $\bullet$ & $\bullet$  &$\otimes$ & &  &&\\
\hline  $\bullet$ & $\bullet$  &$\bullet$  &$\bullet$  & & &\\
 \hline

\end{tabular}
\quad $\Rightarrow$\quad
\begin{tabular}{|p{0.07cm}|p{0.07cm}|p{0.07cm}|p{0.07cm}|p{0.07cm}|p{0.07cm}|p{0.07cm}|}
\hline & & &  & & &   \\
\hline $+$& & & &  & &  \\
\hline $+$&$+$ & & &  & &  \\
\hline $\otimes$ &$-$ &$-$ & & & &  \\
\hline $\bullet$ & $\otimes$    & $-$  &$+$ &  && \\
\hline $\bullet$ & $\bullet$  &$\otimes$ &$\otimes$ &$-$  &&\\
\hline  $\bullet$ & $\bullet$  &$\bullet$  &$\bullet$  & & &\\
 \hline

\end{tabular}
\quad $\Rightarrow$\quad
\begin{tabular}{|p{0.07cm}|p{0.07cm}|p{0.07cm}|p{0.07cm}|p{0.07cm}|p{0.07cm}|p{0.07cm}|}
\hline & & &  & & &   \\
\hline $+$& & & &  & &  \\
\hline $+$&$+$ & & &  & &  \\
\hline $\otimes$ &$-$ &$-$ & & & &  \\
\hline $\bullet$ & $\otimes$    & $-$  & $+$&  && \\
\hline $\bullet$ & $\bullet$  &$\otimes$ &$\otimes$ &$-$  &&\\
\hline  $\bullet$ & $\bullet$  &$\bullet$  &$\bullet$  &$\otimes$ & &\\
 \hline

\end{tabular}}
\end{center}
\vspace{0.3cm}
\begin{center}
\normalsize {\quad $\Rightarrow$\quad $\Dc$ =
\begin{tabular}{|p{0.07cm}|p{0.07cm}|p{0.07cm}|p{0.07cm}|p{0.07cm}|p{0.07cm}|p{0.07cm}|}
\hline & & &  & & &   \\
\hline $+$& & & &  & &  \\
\hline $+$&$+$ & & &  & &  \\
\hline $\otimes$ &$-$ &$-$ & & & &  \\
\hline $\bullet$ & $\otimes$    & $-$  & $+$&  && \\
\hline $\bullet$ & $\bullet$  &$\otimes$ &$\otimes$ &$-$  &&\\
\hline  $\bullet$ & $\bullet$  &$\bullet$  &$\bullet$  &$\otimes$ &$\otimes$ &\\
 \hline

\end{tabular}}

\end{center}

The subset $C(D)$ consists of  $D$ and two roots $(6,4)$ and $(7,6)$.
In this example,  $D'=\{(6,4)\}$ and $C(D) \ne D\cup D'$. The invariant, corresponding to  $(7,6)$ can't be constructed by the method of the paper  \cite{A2}.

To each root of $C(D)$, we attach the polynomial
$F_{4,1} = x_{41}$,~  $F_{5,2} = x_{52}$,~  $F_{6,3} = x_{63}$,~
 $F_{6,4} = x_{64}x_{41}+x_{63}x_{31}$,~
 $F_{7,5} = x_{75}$, ~ $F_{7,6} = x_{76}(x_{64}x_{41}+x_{63}x_{31}) +x_{75}(x_{54}x_{41}+x_{53}x_{31}+x_{52}x_{21})$.

The basic sell $V_D^\circ$ consists of all $X\in \nx_-$ obeying relations
$ x_{\eta}=0$, where $\eta$  is filled with "$\bullet$" on the diagram, and  $F_{4,1} \ne 0$, ~ $F_{5,2} \ne 0$,~  $F_{6,3} \ne 0$,~
$F_{7,5} \ne 0$.

The field of invariants  of the basic cell  $V_D^\circ$
is the field of rational functions of $F_{4,1},~ F_{5,2},~ F_{6,3}, ~F_{6,4}, ~F_{7,5},~ F_{7,6}$.

The field of invariants  of the basic variety
 $V_{D,\vphi}$ is the field  of rational functions of  $F_{6,4}, F_{7,6}$.

\end{document}